\pdfoutput=1

\documentclass[12pt]{elsarticle}
\usepackage[utf8]{inputenc}
\usepackage{graphicx}
\usepackage{amssymb}
\usepackage{tabto}
\usepackage{amsmath}
\usepackage{mathtools}
\usepackage{natbib} 
\usepackage{booktabs}
\usepackage{pgfplotstable}
\usepackage{verbatim}
\usepackage{multirow}
\usepackage{varwidth}
\usepackage{setspace}
\newcommand{\answer}[1]{\textcolor{black}{#1}}

\usepackage[T1]{fontenc}
\usepackage{tikz}
\usetikzlibrary{shapes.geometric,arrows,calc,positioning}

\tikzstyle{startstop} = [rectangle, rounded corners, minimum width=3cm, minimum height=1cm,text centered, draw=black, fill=red!30]
\tikzstyle{io} = [trapezium, trapezium left angle=70, trapezium right angle=110, minimum width=3cm, minimum height=1cm, text centered, draw=black, fill=blue!30]
\tikzstyle{process} = [rectangle, minimum width=3cm, minimum height=1.5cm, text centered, text width=3cm, draw=black, fill=orange!30]
\tikzstyle{decision} = [diamond, minimum width=3cm, minimum height=1cm, text centered, draw=black, fill=green!30]
\tikzstyle{arrow} = [thick,->,>=stealth]

\begin{document}
\begin{frontmatter}

\author[1,2]{Cyril Koch}
\ead{cyril.koch@utt.fr}
\author[1]{Taha Arbaoui} 
\ead{taha.arbaoui@utt.fr}
\author[1]{Yassine Ouazene} 
\ead{yassine.ouazene@utt.fr}
\author[1]{Farouk Yalaoui}
\ead{farouk.yalaoui@utt.fr}
\author[2]{Humbert De Brunier} 
\author[2]{Nicolas Jaunet} 
\author[2]{Antoine De Wulf}

\address[1]{Computer Laboratory and Digital Society (LIST3N), University of Technology of
Troyes, 12 Rue Marie Curie, 10000 Troyes, France}
\address[2]{Manufacture Française des Pneumatiques Michelin, 7 Avenue Président Coty, 10600,
La-Chapelle-Saint-Luc, France}

\title{A Matheuristic Approach for Solving a Simultaneous Lot Sizing and Scheduling Problem with Client Prioritization in Tire Industry}

\begin{abstract}
This paper introduces an integrated lot sizing and scheduling problem inspired from a real-world application in off-the-road tire industry. This problem considers the assignment of different items on parallel machines with complex eligibility constraints within a finite planning horizon. It also considers a large panel of specific constraints such as:  backordering, a limited number of setups, upstream resources saturation and customers prioritization. 
A novel mixed integer formulation is proposed with the objective of optimizing different normalized criteria related to the inventory and service level performance. Based on this mathematical formulation, a problem-based matheuristic method that solves the lot sizing and assignment problems separately is proposed to solve the industrial case.
A computational study and sensitivity analysis are carried out based on real-world data with up to 170 products, 70 unrelated parallel machines and 42 periods. The obtained results show the effectiveness of the proposed approach on improving the company's solution. \answer{Indeed, the two most important KPIs for the management have been optimized of respectively 32\% for the backorders and 13\% for the overstock. Moreover, the computational time have been reduced significantly.}

\end{abstract}

\end{frontmatter}

\section{Introduction}

In today's globalized just-in-time economy, production systems are facing greater challenges. Customers demand more specialized products which leads to increasing products varieties and demand variations. To deal with these challenges, companies need sophisticated decision-support systems. In addition, constraints to consider change significantly from one industry to another. This paper tackles the simultaneous lot sizing and scheduling problem for the tire industry.

Lot sizing problems have been widely studied in the literature over the last decades. The expected output of lot sizing is to give a complete picture of how many pieces to produce  and how many pieces to carry in inventory at each period over a planning horizon. It takes its origin in the well-known Economic Order Quantity (EOQ) model \cite{Harris} under the assumption of single item, constant demand and infinite planning horizon. Since then, numerous researchers have built more realistic models to cope with real-world problems. The production capacity limitation is a significant constraint that production managers have to cope with. The Capacitated Lot Sizing Problem (CLSP) has been proven to be NP-hard by Bitran \cite{bitran1982computational}. Since then, various extensions of the lot sizing problem have been studied extensively and can be classified based on several criteria. An exhaustive review on the CLSP can be found in \cite{drexl1997lot} and \cite{pinedo2005planning} and a classification of criteria is presented in \cite{brahimi2017production}. We limit our study to the multi-item version of the CLSP with a focus on the single-level simultaneous lot sizing and scheduling problem (see \cite{Worbelauer2019Simultaneous} and \cite{copil2017simultaneous} for a more exhaustive classification). Please refer to \cite{brahimi2017single} for single-item LSP and to \cite{buschkuhl2010dynamic} for multi-level LSP.

Our contribution is threefold. First, we tackle an industrial case and demonstrate the efficiency of our method to solve it. Second, twenty specific constraints such as the number of setup per period, upstream resources saturation and customer prioritization are presented. To the best of our knowledge, these constraints have never been modeled and considered together in the literature. Finally, we propose a MIP formulation and a problem-based matheuristic to solve the problem efficiently. Several objectives are taken in consideration and objective function parameters are calibrated thanks to the Taguchi procedure. A sensitivity analysis of one particular parameter is also conducted.

The remainder of this paper is organized as follows. Section \ref{Literature_review)} provides a literature review of the considered problem. The problem description and the assumptions of this study are described in Section \ref{Problem description}. Section \ref{MIP formulation} describes the mathematical formulation proposed to model the planning problem. In Section \ref{Matheuristic}, a matheuristic method based on an original hybrid sequential approach is introduced to solve the problem. A model size comparison is also presented. Section \ref{Computational Results} details the real case study results based on large performance testing campaigns and sensitivity analysis. Finally, Section \ref{conclusion} summarizes the study and gives future research directions.

\section{Literature Review}
\label{Literature_review)}
As mentioned before, the underlying constraints to consider vary significantly from one industry to another. Recent papers dealing with simultaneous lot sizing and scheduling (LSS) problems put more focus on particular features of these industries \cite{copil2017simultaneous}: the beverage industry \cite{toledo2015synchronized}, \cite{toledo2015relax} and \cite{baldo2017alternative}, steel manufacturing \cite{li2017fix}, \cite{de2008lot} and \cite{de2007joint}, the automotive industry  \cite{deeratanasrikul2017multiple} and \cite{gnoni2003production}, glass-container manufacturing \cite{toledo2016mathematical}, the tile industry  \cite{ramezanian2017simultaneous}, the chemical industry \cite{cunha2018integrated}, and the paper industry \cite{figueira2013hybrid}. Additionally some research has been conducted in the tire industry \cite{lasdon1971efficient} and \cite{jans2004industrial} and more specifically in the off-the-road tire industry \cite{degraeve1997tire} and \cite{koch2020dedicated}.

Some researchers used meta-heuristic to deal with the studied NP-hard problem. One of the most popular meta-heuristic is the genetic algorithm (GA). Babaei et al. \cite{babaei2014genetic} applied it on the LSS problem in capacitated flow shop environment with consideration of backlogging and sequence-dependent setups. They used a heuristic proposed by Mohammadi et al. \cite{mohammadi2011genetic} to set the initial population and crossover operators from Ruis and Maroto \cite{ruiz2005solving}. Vincent et al. \cite{vincent2020population} recently presented a population-based heuristic to solve a multi-item capacitated lot sizing problem with setup time and unrelated parallel machines. A dynamic constructive heuristic was proposed to generate a set of initial solutions. A diversification procedure was realised using a path-relinking strategy and finally the population was intensified with a local search method. Another well known meta-heuristic is the Simulated Annealing (SA) algorithm. Ceschia et al. \cite{ceschia2017solving} introduced a SA approach to cope with the multi-item single-machine LSS problem. The search space was a vector $V$ (size: number of periods) with values $v_t$ representing the items to produce. The neighbourhood relation considered was the composite Swap $\bigcup$ Insert firstly proposed by Della Croce \cite{della1995generalized}. They also provided a hybrid method combining a MILP formulation and SA. The SA was run for a short time period and then the solution was injected as initial solution for the MILP using \textit{warm start} functionality of the CPLEX solver.

Beyond general meta-heuristic approaches, an increasing number of researchers bring forward decomposition heuristics. They usually start from a new MIP formulation from a real-world case study and propose approaches that can be product-, machine- or time based- decomposition heuristics. For instance, Meyr and Mann \cite{meyr2013decomposition} described a heuristic for the general LSS problem for parallel production lines with consideration of backlogging. The multi-line problem was divided into a series of single-line problems easier to solve. They discussed priority rules for product assignment to avoid setups and save production and inventory costs. Time-based decomposition approaches are much more widespread, mostly with Relax-and-Fix (R\&F) and Fix-and-optimize (F\&O) heuristics. De Araujo et al. \cite{de2008lot} applied a R\&F approach to solve an integrated two-level lot sizing and furnace scheduling problem in small foundries. They presented a MIP formulation for the case study and solved each macro-period divided in 10 micro periods one after another. A descent heuristic procedure was also proposed and improved using diminished neighbourhood search and simulating annealing. Rodoplu et al. \cite{rodoplu2020fix} used a R\&F heuristic to address a single-item lot sizing problem with a flow shop system and energy constraints. A worst-case analysis of R\&F algorithm for lot sizing problems have recently been conducted by Absi and Van den Heuvel \cite{absi2019worst}. They analysed the impact of time window length, effect of overlapping time windows and capacity constraints and showed that even for simple instances with time-invariant parameters, the worst-case ratio may be unbounded. R\&F procedures are often coupled with F\&O heuristics as in \cite{deeratanasrikul2017multiple} where a multi-level multi-machine CLSP with setup times was dealt with. The R\&F heuristic proposed was a two-level partition of the problem. Each sub-problem related to a time-window of $\lambda$ consecutive periods in one stage. The R\&F solution was used as an initial solution for the F\&O procedure which was based on machine-decomposition strategy. The first to introduce F\&O were Helbert and Sahling \cite{helber2010fix} and they also applied it to the multi-level CLSP. Toscano et al. \cite{toscano2020formulation} have recently developed a F\&O algorithm to solve a synchronized two-stage CLSP with mandatory temporal cleaning and sequence-dependent changeover in a soft drink company. F\&O heuristics have also been combined with local search such variable neighborhood search as in \cite{chen2015fix} and \cite{li2017fix}. Li et al. \cite{li2017fix} dealt with a multi-item lot sizing problem in the steel industry with demand class formulation and stochastic demand where backlogging and overtime costs were also incurred. An original F\&O approach based on ``k-degree decomposition'' that proposes several decomposition among  products, resources and time horizon was presented. Then they proposed an integrative F\&O and VNS procedure. The demand class formulation is an original contribution of this paper. The work of Gruson et al. \cite{gruson2018impact} is also to be noted on this topic. They provided a service level analysis for lot sizing problems with backlogging. They introduced two service-level definitions and modelling based on backorder fix and variable costs. Moreover, Gören and Tunali \cite{goren2015solving} used a sequential hybrid approach with F\&O and GA for the multi-level CLSP with setup carryover. They used a time window of five consecutive periods and the sub-problems were then solved with a GA. Besides, Toledo et al. \cite{toledo2015relax} and Baldo et al. \cite{baldo2017alternative} both brought forward for the soft drink and brewery industry MIP-based heuristics with binary variables relaxation for two-level multi-item and multi-machine CLSP with sequence-dependent setup at each level. In \cite{toledo2015relax} they suggested two relaxation of binary variables into continuous variables to minimize setup, production and inventory costs. In \cite{baldo2017alternative}, they added backorder consideration and coped with two-level problem with tanks providing bottling lines. The MIP-based sequential approach focused on the decomposition of stages to minimize setup, production, inventory and backorder costs. Finally, an interesting fuzzy mathematical programming and self-adaptating artificial fish swarm algorithm for just-in-time energy aware flow shop scheduling problem with outsourcing option was presented by Tirkolaee et al. \cite{tirkolaee2020fuzzy}. They proposed a novel bi-objective mixed-integer linear programming model (MILP). They treated it as a single objective MILP using a multi-objective fuzzy mathematical programming technique and implemented a sensitivity analysis to study the behavior of the objectives with real-world conditions.

\section{Problem Description}
\label{Problem description}
\subsection{Industrial Process}
\answer{
This work tackles a planning problem in tire manufacturing. A tire is a complex structure composed of many layers. The main components can be summarized as follow (this list is non-exhaustive, see Figure 1): (i) The inner liner plays the role of an inner tube; (ii) One to three casing plies. It is also called the carcass and refers to the main body of the tire; (iii) The crown ply. This belt is placed on top of the casing and consists of layers of fabric to avoid tire deformation; (iv) The tire tread. It is the part of the tire that comes into direct contact with the road; (v) The sidewalls are found at the sides of the tires. It protects the carcass and prevents tears or punctures on the side of the tire; (vi) The tire bead is composed of steel wire wrapped into a thin layer of rubber. These large steel cords are wound together to form a cable to compose the bead cores.}

\begin{figure}[!h]
	\centering
	\includegraphics[scale=0.7]{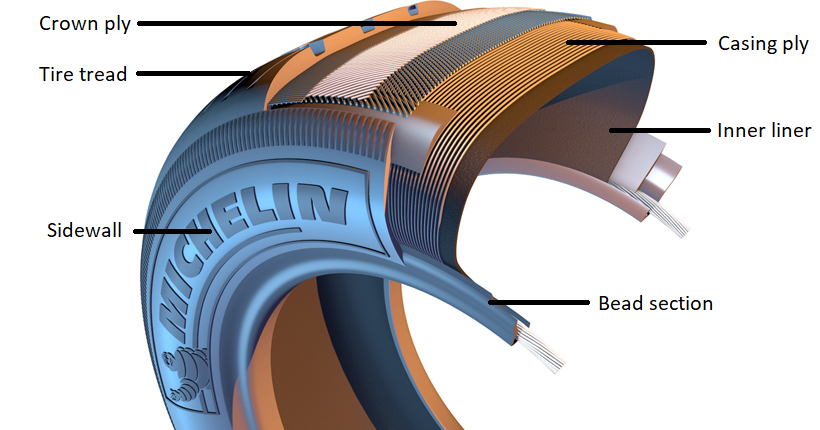}
	\caption{Tire architecture}
\end{figure}

\answer{
The whole production process can be divided in five major sub-processes (see Figure 2\footnote{Tire architecture description provided by the company}) \cite{boudha_kumar_srimannarayana_koyan_2011}. The first subprocess is the production of a homogeneous rubber material based on  four main compounds: elastomers, renforcing fillers, platicizers and other chemicals elements. These components are mixed in a banbury mixer to obtain a homogeneous rubber material in the form of thin layers. In addition to the rubber compounds, two more raw materials are necessary to build a tire: textile reinforcements and steel wires, ensuring its rigidity and geometry. The second sub-process is the production of semi-finite products using two different technologies: first the extruding and calendaring process; and second the profiling and cutting process. Then the task of building the tire begins. This third sub-process is called the assembling sub-process. The tire building machine also needs a resource called a drum (a rotative cylinder) on which the different parts of the tire are assembled. The inner liner, the casing plies,the bead cores and the sidewalls are put together. After that the building machine operates the conformation operation. The edges of the drum are brought together and the center is inflated: this operation gives the final toroidal shape of the tire. Finally the crown ply and the tire tread are incorporated in that structure to obtain what is called a ``green'' or uncured tire. Fourth is the curing and vulcanizing process. A tire-specific mold is placed into the curing press, which is also called a heater or curing press, and the green tire is put into that mold. A bladder filled with pressurized hot fluid in the center of the mold forces the still malleable substance of the green tire to flow into all the cavities of the tread pattern engraved inside the mold. The heat of the fluid starts the curing process. The increase of temperature causes the sulfur contained in the rubber compound to bound with the rubber molecules. This is what we call vulcanization. The rubber is then transformed from a plastic to an elastic state. When ejected from the mold and after cooling the tire has taken on his final shape and properties. Finally inspection and finishing operations remains before the tire is stored in the warehouse.
}
\begin{figure}[!h]
	\centering
	\includegraphics[scale=0.8]{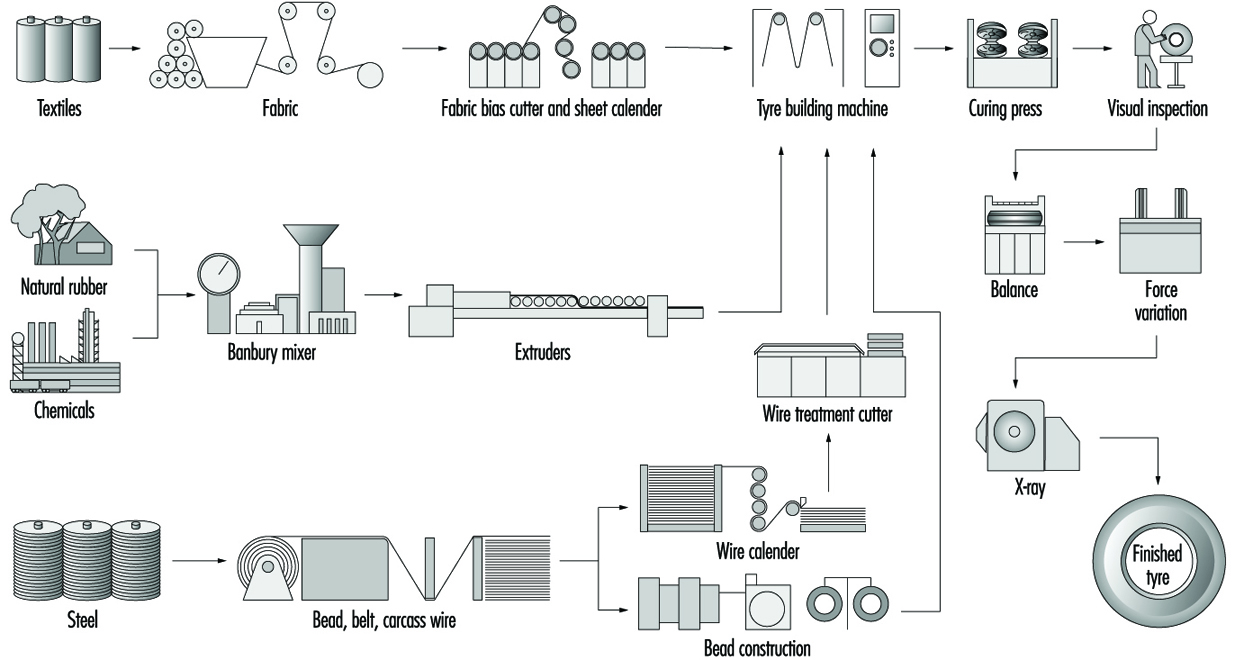}
	\caption{Tire Manufacturing Process from \cite{boudha_kumar_srimannarayana_koyan_2011}}
\end{figure}

\subsection{Curing Workshop Planning}
\answer{To deal with the planning problem of this plant, the company should adress the balance between the assembling and curing subprocesses capacity. Indeed this production environnement can be described as a hybrid flow shop environnement with two stages and a different number of machines at each stage with tire-dependent processing times with high variability. The assembling workshop has twenty-two assembling machines for a mean processing time around fourty minutes while seventy curing presses are available for a mean processing time of two hours. The throughput analysis of the production process provided by the company allows to identify the curing workshop as the bottleneck. However, depending on the mix of tires to produce, the bottleneck can move from curing to assembling workshop from time to time. Therefore, the stability of the production process between assembling and curing workshop is balanced on a knife edge. Hence, the idea is to organise the production of the plant based on the curing workshop planning problem while considering at the same time the saturation of assembling workshop.}\\

\answer{
We therefore focus on the curing workshop simultaneous planning and scheduling problem. There is a wide portfolio of around 170 tires to produce on 70 curing presses over 42 periods. In addition, the portfolio is getting wider to match customer expectations and makes the considered problem more difficult. The production is based on a make-to-stock inventory policy, so that the inventory level stays between a minimum and a maximum level calculated to prevent shortage and keep Working Capital Requirement to a minimum. The number of campaign endings per week is limited to avoid too much raw material loss in semi-finite upstream workshop. Also, the number of setups per period is limited by human ressources. A new mold setup requires consequent setup times (20\% of the usual daily production yield) and the planning problem is highly restricted by the eligibility matrix between tires and curing presses.\\}

During the curing process, the green tire is put into a mold that provides a specific pattern for the tire. Each mold is tire-specific: it can be used for exactly one type of tire. For some tire references several molds are available, though for most tires there is only one mold. Every mold can be placed in several curing presss, respecting the eligibility matrix. Nonetheless, each press can contain at most one mold at a time. The curing time depends on the tire produced and the curing press used. The curing presses capacity therefore links together different tire references that compete for the same resource - available time of a given curing press where the molds can be placed in. Except for the first and the last period of the production campaign, tires are produced in a continuous run and production is always done at full capacity. This type of production is often referred to as ``all-or-nothing'' production. Also, only one type of tires can be cured in a curing press within one period. Thus, our problem is classified as a small-bucket lot sizing problem.

\section{MIP formulation}
\label{MIP formulation}
The production planning problems encountered in the industry may be intractable in numerous situations due to several practical constraints. The demand over the planning horizon is known in advance (deterministic). In order to deal with situations when demand cannot be met in time, the company allows backlogging. Specific constraints are also added such as the number of ``campaign endings'' within the planning horizon or the number of different tires produced at each period. In Figure 3 is presented in the red circles the constraints to be considered. Blue circles denote the different sources where the constraints comes from. The objective of the proposed approach is to find a good feasible solution for the single-level multi-item multi-machine with deterministic demand, backlogging, sequence independent setup times and specific constraints problem. The method also takes in consideration different classes of demand representing client prioritization.
\begin{figure}[!h]
	\centering
	\includegraphics[scale=1]{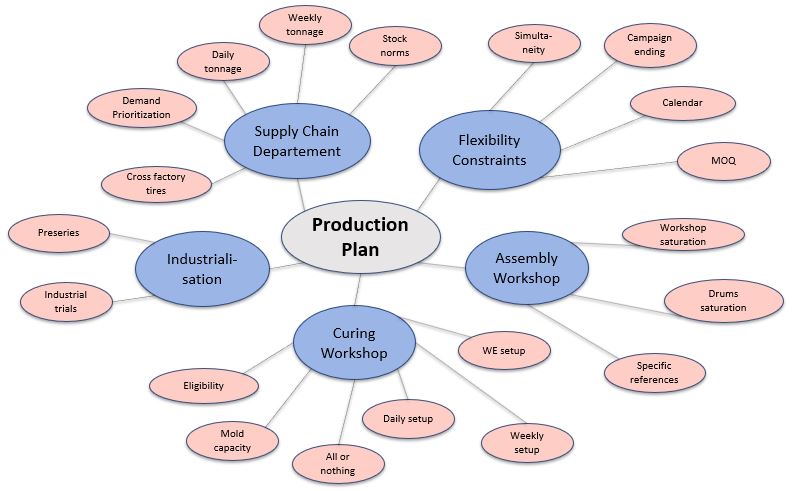}
	\caption{Constraints to be considered and variety of sources diagram}
\end{figure}
 
\textbf{Indexes and Sets}\\
$A$: Number of tires in the portfolio, $a$ = 1~..~$A$\\
$N$: Number of items to plan, $i$ = 1~..~$N$, $N\ge A$\\
$N_a$: Number of items $i$ for tire $a$, $i$ = 1$~$..$~N_a$\\
$P$: Number of curing presses, $p$ = 1~..~$P$\\
$T$: Number of micro-periods $t$ in one macro-period $h$, $t$ = 1~..~$T$\\
$H$: Number of macro-period, $h$ = 1~..~$H$\\
$W$: Number of workshops in assembling shop, $w$ = 1~..~$W$\\
$N_w$: Set of items $i$ doable on workshop $w$\\
$N_d$: Number of assembling machine resource types (drums), $d$ = 1~..~$N_d$\\
$C$: Set of demand classes for client prioritization for tire $a$ in 1..$A$, $c$= ${ C_1..C_\gamma}$\\

\textbf{Parameters}\\
Inventory and capacity parameters\\
$\overline{S_{at}}$: Overstock of tire $a$ at period $t$\\
$\underline{S_{at}}$: Understock of tire $a$ at period $t$\\
$D_{act}$: Demand of class $c$ for tire $a$ at period t\\
$\gamma_{c}$: Weight vector to balance demand class priorities\\
$K_a$: Number of curing press resources (molds) available for tire $a$\\
$(M_{ap})$: Eligibility matrix tire - curing press\\
$R_{it}$: Daily rate of a curing press for item $i$ during period $t$\\

Weights parameters \\
$V_{t}$: Targeted weight to produce within one micro-period $t$\\
\answer{
$\overline{v_t}$: Upper tolerance for targeted weight at micro-period $t$\\
$\underline{v_t}$: Lower tolerance for targeted weight at micro-period $t$\\
$\overline{v_h}$: Upper tolerance for targeted weight at macro-period $h$\\
$\underline{v_h}$: Lower tolerance for targeted weight at macro-period $h$\\
}
$\Omega_i$: Unit weight of item $i$\\

Upstream workshop saturation parameters\\
$\overline{S_w}$: Maximum saturation of workshop $w$ in assembling shop (time unit)\\
$N_{id}$: Set of item $i$ that can be produced with drum $d$ in assembling shop \\
$T_i$: Unit time of production of item $i$ in assembling shop \\
$K_d$: Number of assembling shop resources (drums) of type $d$ \\
$\varepsilon_{id}$: Number of molds for a given item $i$ that can be furnished by the production of one drum $d$ in upstream assembling shop\\

Flexibility parameters\\
$S$: Maximum number of items $i$ produced simultaneously per micro-period $t$\\
\answer{
$\overline{C_t}$: Maximum number of molds setup per micro-period $t$\\
$\overline{C_h}$: Maximum number of molds setup per macro-period $h$\\
$\overline{e_h}$: Maximum number of campaign ending per macro-period $h$\\ 
$\tau_m$: Minimum number of micro-periods $t$ of production after a setup \\
$\tau_s$: Duration of campaign suspension to count a mold setup\\
$\tau_e$: Duration of campaign suspension to count a campaign ending\\
}

Problem specific parameters\\
$(U_{pt})$: Matrix of unavailable curing press $p$ for maintenance in period $t$\\
$(P_{ipt})$: Matrix of production of item $i$ on curing press $p$ in micro-period $t$ to be enforced for industrial trial \\
$(F_{t})$: Vector of days off $t$ for each macro-period $h$\\
$\underline{Molds}$ : Minimum number of molds to plan at the same time\\
$Spec$: Set of specific items $i$ that needs to be produced with at least $\underline{Molds}$ molds at a time for quality reasons\\
$M = \max(R_{it})$: A large number that major the maximum yield $R_{it}$\\
\\
\textbf{Decision Variables}\\
$I_{at}$: Inventory level for tire $a$ at the end of micro-period $t$\\
$B_{act}$: Backorder level of class $c$ for tire $a$ at the end of micro-period $t$\\
$X_{ipt}$: Quantity of item $i$ to produce on curing press $p$ in micro-period $t$\\
\answer{$\Delta_{idt}$: Number of drums of type $d$ needed for item $i$ in micro-period $t$} \\
$Y_{ipt}$: Binary variable that equals 1 if there is a new production of an item $i$ on curing press $p$ at period $t$; 0 otherwise\\
\answer{$\sigma_{it}$: Binary variable that equals 1 if item $i$ is being cured at micro-period $t$; 0 otherwise\\
$s_{ipt}$: Binary variable that equals 1 if there is a mold setup for item $i$ on curing press $p$ in micro-period $t$; 0 otherwise\\
$e_{it}$: Binary variable that equals 1 if there is a campaign ending for item $i$ in micro-period $t$; 0 otherwise\\
$m_{ipt}$: Binary variable that equals 1 if in micro-period $t$ the minimum number of days of production after a setup have not been reached on curing press $p$ ; 0 otherwise\\
$\delta_{idt}$: Binary variable that equals 1 if in micro-period $t$ the drum $d$ is used to produce item $i$; 0 otherwise\\
}
\\
The objective of the production planner of the tire company is to optimize several criteria. The main objective is to prevent shortage. We consider $n$ demand classes ${ C_1}..{ C_\gamma}$. All demand of class ${ C_1}$ should be satisfied before demand of class ${ C_2}..{ C_\gamma}$, and so on. All demand of class ${ C_2}$ should be satisfied before demand of class ${ C_3}..{ C_\gamma}$. Once backordering has been minimized, the production planner tries to keep every tire between a minimum and a maximum inventory level set by the supply chain department. $\lambda_j$ coefficients allows the right balancing between chosen KPIs. The objective function is expressed as follows.\\\\
Minimize:
\begin{equation} \label{1}
   Z = \sum_{c \in C}\lambda_{c}\sum_{a=1}^{A}\sum_{t=1}^{T}B_{act} ~+~ \lambda_{\overline{S}} \sum_{a=1}^{A}\sum_{t=1}^{T}\overline{S}_{at} ~+~ \lambda_{\underline{S}} \sum_{a=1}^{A}\sum_{t=1}^{T}\underline{S}_{at}   
\end{equation}\\
To avoid scaling effects, another objective function with normalization coefficients $\mu_j$ is formulated below.\\

Minimize:
\begin{equation}\label{2}
Z'= \sum_{c \in C}\frac{\lambda_{c}}{\mu_c} \sum_{a=1}^{A}\sum_{t=1}^{T}B_{act}
~+~ \frac{\lambda_{\overline{S}}}{\mu_{\overline{S}}} \sum_{a=1}^{A}\sum_{t=1}^{T}\overline{S}_{at}
~+~  \frac{\lambda_{\underline{S}}}{\mu_{\underline{S}}} \sum_{a=1}^{A}\sum_{t=1}^{T}\underline{S}_{at} 
\end{equation}\\

Where
\begin{equation}\label{3}
 \mu_c=\sum_{a=1}^{A}B_{ac0}+\sum_{t=1}^{T}D_{act} , c \in C  
\end{equation}
\begin{equation}\label{4}
\mu_{\overline{S}}=\sum_{a=1}^{A}(\overline{S}_{a0}+K_a R_{a} T)
\end{equation}
\begin{equation}\label{5}
\mu_{\underline{S}}=\sum_{a=1}^{A}\underline{S}_{a0} +  \sum_{c \in C}\mu_{c}    
\end{equation}

Subject to
\begin{equation}\label{6}
\begin{split}
I_{at-1} + \sum_{i=1}^{N_a}\sum_{p=1}^{P} X_{ipt} = I_{at} + \sum_{c \in C} (D_{act} - B_{act} + B_{act-1}) , a = 1..A, t = 1..T
\end{split}
\end{equation}
\begin{equation}\label{7}
I_{at-1}  + \sum_{i=1}^{N_a}\sum_{p=1}^{P} X_{ipt} \ge I_{at} , a = 1..A, t = 1..T
\end{equation}
\begin{equation}\label{8}
I_{at-1} + \sum_{i=1}^{N_a}\sum_{p=1}^{P} X_{ipt} \ge D_{act} - B_{act} + B_{act-1} , a = 1..A, c = C_1 .. C_{n-1}, t = 1..T
\end{equation}
\begin{equation}\label{9}
B_{act} \le D_{act} + B_{act-1} , a = 1..A, c = C_2 .. C_{n}, t = 1..T
\end{equation}\\
Constraint (\ref{6}) is the inventory balance equation, with consideration of backorder. Thanks to Constraints (\ref{7}) to (\ref{9}) demand prioritization is correctly handled. Constraint (\ref{7}) avoid the creation of ``ghost''  inventory through inventory Equation (\ref{6}). Constraint (\ref{8}) ensures the propagation of backorder prioritization from one period to another in inventory Constraint (\ref{6}). Finally Constraints (\ref{9}) avoid the creation of ``ghost''  backorder through inventory Equation (\ref{6}).\\
\begin{equation}\label{10}
X_{ipt}\le M Y_{ipt}, i = 1..N, p = 1..P, t = 1..T
\end{equation}
\begin{equation}\label{11}
X_{ipt}\ge Y_{ipt}, i = 1..N, p = 1..P, t = 1..T
\end{equation}
\begin{equation}\label{12}
\sum_{i=1}^{N} Y_{ipt}\le 1, p = 1..P, t = 1..T
\end{equation}
\begin{equation}\label{13}
\sum_{i=1}^{N_a}\sum_{p=1}^P Y_{ipt}\le K_a, a = 1..A, t = 1..T
\end{equation}
\begin{equation}\label{14}
X_{ipt}= R_{it}  Y_{ipt} , i = 1..N, p = 1..P, t = 1..T
\end{equation}
Constraints (\ref{10}) and (\ref{11}) link the setup binary variable and the production variable together. 
Constraint (\ref{12}) allows only one item to be produced on a curing press within a single period. 
Constraints (\ref{13}) and (\ref{14}) represent capacity constraints. 
Given a particular tire $a$, Constraint (\ref{13}) limits the number of curing presses used to the mold capacity of that tire ($K_a$).
Constraint (\ref{14}) combined with (\ref{12}) describes the ``all-or-nothing'' policy. 
The production rate $R_{it}$ can be affected by events planned in advance in the calendar such operator training or meetings imposed by the management.\\
\begin{equation}\label{15}
s_{ipt}\ge Y_{ipt} - \sum_{o=t-\tau_s}^{t-1} Y_{ipo} , i = 1..N, p = 1..P, t = 1..T
\end{equation}
\begin{equation}\label{16}
\begin{split}
F_{t-1}   (Y_{ipt}-Y_{ipt-1-\sum_{o =t-\tau_m}^{t-1}F_{o}} -1 + \frac{1}{\tau_s N} \sum_{o=t-\tau_s}^{t-1}\sum_{j\neq i}Y_{jpo}) \\ \le F_{t-1}   s_{ipt} ,i = 1..N, p = 1..P, t = 1..T
\end{split}
\end{equation}
\begin{equation}\label{17}
\begin{split}
 (1-F_{t-1})   (Y_{ipt}-Y_{ipt-1} -1 + \frac{1}{\tau_s N} \sum_{o=t-\tau_s}^{t-1}\sum_{j\neq i}Y_{jpo}) \\ \le  (1-F_{t-1})   s_{ipt} ,i = 1..N, p = 1..P, t = 1..T
\end{split}
\end{equation}

\begin{equation}\label{18}
\sum_{i=1}^{N}\sum_{p=1}^P s_{ipt}\le \overline{C_t}, t = 1..T
\end{equation}
\begin{equation}\label{19}
\sum_{i=1}^{N}\sum_{p=1}^P\sum_{t=1}^{T} s_{ipt}\le \overline{C_h}
\end{equation}\\

Each item changeover from one period to another incurs a mold setup, which requires a human resource and affects the curing press productivity during the changeover. To respect the changeover capacity two parameters have been set: the maximum number of mold setup per day ($\overline{C_t}$) and per week ($\overline{C_h}$). Constraints (\ref{15}) to (\ref{18}) limit the number of mold setup per period (days) whereas Constraints (\ref{15}), (\ref{16}), (\ref{17}) and (\ref{19}) limit the number of mold setup per set of periods (weeks). In Constraint (\ref{15}), a mold setup is incurred if no production have been done for a given item $i$ on curing press $p$ within the last $\tau_s$ periods. This check over $\tau_s$ periods in the past is necessary as a curing press is allowed not to produce an item during $\tau_s$ periods without considering a mold setup if the mold stays in the curing press (see Figure 4). It means that the curing press cannot be used to produce another item during these periods. Constraint (\ref{15}) forces to count a mold setup if during a campaign suspension this time slot have been used to produce another (short) campaign. In Constraint (\ref{16}) and (\ref{17}) a parameter $F_t$ allows the consideration of days-off in the calendar.\\

\begin{figure}[!h]
	\centering
	\includegraphics[scale=0.7]{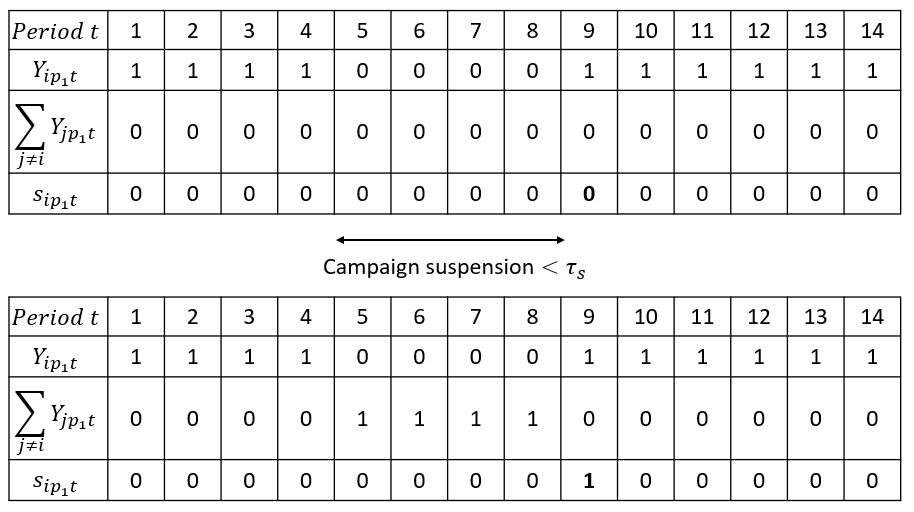}
	\caption{Mold setup illustration}
\end{figure}

\begin{equation}\label{20}
\sum_{i=1}^{N}\sum_{p=1}^P\sum_{t=1}^{T} X_{ipt} \Omega_i\le \sum_{t=1}^{T}V_t + \overline{v_h}
\end{equation}
\begin{equation}\label{21}
\sum_{i=1}^{N}\sum_{p=1}^P\sum_{t=1}^{T} X_{ipt} \Omega_i\ge \sum_{t=1}^{T} V_t - \underline{v_h}
\end{equation}
\begin{equation}\label{22}
\sum_{i=1}^{N}\sum_{p=1}^P X_{ipt} \Omega_i\le V_t + \overline{v_t}, t = 1..T
\end{equation}
\begin{equation}\label{23}
\sum_{i=1}^{N}\sum_{p=1}^P X_{ipt} \Omega_i\ge V_t - \underline{v_t}, t = 1..T
\end{equation}

Every week, the factory commits to the supply chain department a targeted total weight of items to produce. 
The workforce needed is determined to satisfy this target. 
Respecting these constraints is very important to keep the cost of production the lowest possible.
Constraints (\ref{20}) to (\ref{23}) ensure the targeted total weight of items produced stays between upper and lower bounds, per days and per week, where $\Omega_i$ represents the unit weight of item $i$, $\underline{v_h}$ and $\overline{v_h}$ weekly bounds and $\underline{v_t}$ and $\overline{v_t}$ daily bounds.\\

\begin{equation}\label{24}
\sigma_{it} \ge \frac{1}{P}   \sum_{p=1}^P Y_{ipt}, i = 1..N, t = 1..T
\end{equation}
\begin{equation}\label{25}
\sigma_{it} \le  \sum_{p=1}^P Y_{ipt}, i = 1..N, t = 1..T 
\end{equation}
\begin{equation}\label{26}
\sum_{i=1}^N\sigma_{it} \le S, t = 1..T 
\end{equation}

Constraints (\ref{24}) to (\ref{26}) ensure that the number of different items $i$ produced in the same period do not exceed the limit set by the company ($S$). 
This limit is empirical and makes sure that the upstream workshops (assembling and semi-finished products) are not saturated.
\begin{equation}\label{27}
e_{it} \ge \frac{1}{P} \sum_{p=1}^P Y_{ipt-\tau_e} - P  \sum_{o=t-\tau_e+1}^{t}\sum_{p=1}^P Y_{ipo} , i = 1..N, t = 1..T
\end{equation}
\begin{equation}\label{28}
\sum_{i=1}^{N}\sum_{t=1}^{T}e_{it} \le \overline{e_h}, h = 1..H
\end{equation}

When a production campaign comes to an end, all the semi-finished products need to be produced at the exact quantity to minimize material loss. It requires a particular follow up by the agents in the factory and is quite difficult to manage properly. 
Thus, a limit of campaign ending per week $\overline{e_h}$ has been added.
Similarly to the mold setup constraints, when a production campaign is suspended for a few days (maximum $\tau_e$), it does not count as a campaign ending (see Figure 5). \\
\begin{figure}[!h]
	\centering
	\includegraphics[scale=0.7]{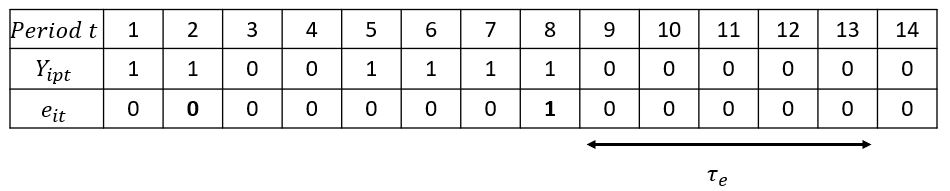}
	\caption{Campaign ending illustration}
\end{figure}

We can also highlight the difference between a production ending on a given press $p$ and a campaign ending as illustrated in Figure 6. Constraints (\ref{27}) and (\ref{28}) state that the number of campaign ending does not exceed a certain limit $\overline{e_h}$ per week.
\begin{figure}[!h]
	\centering
	\includegraphics[scale=0.7]{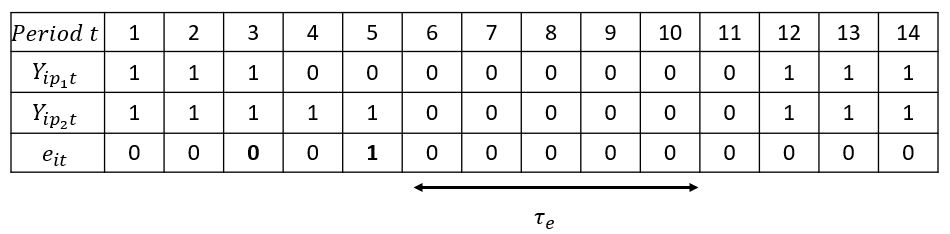}
	\caption{Production VS Campaign ending illustration}
\end{figure}

\begin{equation}\label{29}
\sum_{i\in N_w}\sum_{p=1}^P\sum_{t=1}^{T} X_{ipt} T_i \le \overline{S_w} , w = 1..W
\end{equation}

Thanks to Constraint (\ref{29}), the number of minutes planned on each workshop of the assembling shop does not exceed the saturation limit $\overline{S_w}$. 
This constraint is the reason why the differentiation between two sets of tires $a$ in $1..A$ and items $i$ in $1..N$ is necessary. The same tire $a$ from client view can be made on two different assembling shop workshops. Thus, we need to consider the same tire $a$ as two different items $i$ depending on the shop it is made on. It allows the company to balance the saturation of the different assembling shops and to provide more efficiently the curing workshop bottleneck.\\
\begin{equation}\label{30}
\sum_{i \in N_{id}} \Delta_{idt}\le K_d , d = 1..N_d, t = 1..T
\end{equation}\\
\begin{equation}\label{31}
\Delta_{idt} \ge 2 - M   (1-\delta_{idt}),i = 1..N, d = 1..N_d, t = 1..T
\end{equation}
\begin{equation}\label{32}
\Delta_{idt} \le 2 + M   (1-\delta_{idt}),i = 1..N, d = 1..N_d, t = 1..T
\end{equation}
\begin{equation}\label{33}
\Delta_{idt} \ge - M   \delta_{idt},i = 1..N, d = 1..N_d, t = 1..T
\end{equation}
\begin{equation}\label{34}
\Delta_{idt} \le 1 + M   \delta_{idt},i = 1..N, d = 1..N_d, t = 1..T
\end{equation}

\begin{equation}\label{35}
\Delta_{idt} \ge \sum_{p=1}^P \frac{ Y_{ipt}}{\varepsilon_{id}} ,i = 1..N, d = 1..N_d, t = 1..T
\end{equation}

\begin{equation}\label{36}
 \frac{ \sum_{p=1}^P Y_{ipt}}{\varepsilon_{id}} \ge 1 - M (1 - \delta_{idt}) ,i = 1..N, d = 1..N_d, t = 1..T
\end{equation}
\begin{equation}\label{37}
 \frac{ \sum_{p=1}^P Y_{ipt}}{\varepsilon_{id}} \le 1 + M \delta_{idt} ,i = 1..N, d = 1..N_d, t = 1..T
\end{equation}

There is a limited number $K_d$ of each type of drum $d$ in the assembling shop. Constraint (\ref{30}) ensures that for each drum $d$, the overall production planned in curing process can be absorbed by the assembling shop, where $\Delta_{idt}$ represents the number of drums $d$ necessary for the production of item $i$ planned in the curing workshop in period $t$. Constraints (\ref{31}) to (\ref{37}) help to define $\Delta_{idt}$ variable that can take three values: 0, 1 or 2. These constraints are linearisation constraints.

\begin{equation}\label{38}
(1-F_{t-1})   m_{ipt} \ge (1-F_{t-1})    (Y_{ipt}-Y_{ipt-1}) , i = 1..N, p = 1..P, t = -\tau_m+1 ..T
\end{equation}

\begin{equation}\label{39}
F_{t-1}   m_{ipt} \ge F_{t-1}   Y_{ipt}-Y_{ipt-1-\sum_{o =t-\tau_m}^{t-1}F_{o}} , i = 1..N, p = 1..P, t = -\tau_m+1 ..T
\end{equation}

\begin{equation}\label{40}
(1-F_{t}) Y_{ipt} \ge (1-F_{t})   \sum_{o = t-\tau_m-\sum_{q =t-\tau_m}^{t-1}F_{q}}^{t-1} m_{ipo} , i = 1..N, p = 1..P, t = 1..T
\end{equation}

Constraints (\ref{38}) to (\ref{40}) ensure that once a mold is set up in a given curing press $p$, the production last at least $\tau_m$ days. Again, $F_t$ parameter allows the consideration of days-off in the calendar.\\

\begin{equation}\label{41}
Y_{ipt} \ge P_{ipt} , i = 1..N, p = 1..P, t = 1..T
\end{equation}

\begin{equation}\label{42}
X_{ipt} \le U_{pt}   M, i = 1..N, p = 1..P, t = 1..T
\end{equation}

\begin{equation}\label{43}
\sum_{p=1}^P Y_{ipt} \ge \underline{Mold} ,i \in Spec, t = 1..T
\end{equation}

\begin{equation}\label{44}
X_{ipt} \le M_{ap}   M, a=1..A, i = 1..N_a, p = 1..P, t = 1..T
\end{equation}

The planner regularly needs to enforce some production campaigns trials of a given item $i$ on a curing press $p$ for quality and industrialisation departments. Constraint (\ref{41}) integrates these requirements. On the contrary, Constraint (\ref{42}) enforces that no item is allocated on a curing press shut down for maintenance, industrial trial or any other planned event. Besides when some items are produced (i in $Spec$, set of specific items), it has to be on at least $\underline{Mold}$ molds at a time to avoid quality issues (Constraint (\ref{43})). Constraint (\ref{44}) states that the allocation of items must respect the eligibility matrix.

\begin{equation}\label{45}
Y_{ipt},\sigma_{it}, s_{ipt}, e_{it},\delta_{idt},m_{ipt} \in \{0;1\}
\end{equation}
\begin{equation}\label{46}
X_{ipt},I_{at},B_{act},\Delta_{idt} \in \mathbf{N}
\end{equation}

Finally, Constraint (\ref{45}) defines the boolean decision variables and Constraint (\ref{46}) the production, inventory and number of drums decision variables. To the best of our knowledge, all these constraints have never been modeled and considered together in the literature.

\section{Matheuristic Approach}
\label{Matheuristic}
\answer{When the MIP is applied on  real-world instances, it faces considerable difficulty to solve the problem. To cope with this issue, a problem-specific decomposition-based matheuristic approach is introduced. We decompose our problem into two sub-problems: a lot sizing problem and a machine assignment problem. The aim of this sequential approach is to provide a two-stage algorithm to solve our simultaneous lot sizing and scheduling problem. We firstly find a solution for the lot sizing sub-problem. It consists of putting away the machine assignment decision, and all the associated constraints, of an item $i$ planned at a period $t$. This sub-problem will further be called LSSP (Lot-Sizing Sub-Problem). This relaxation simplifies considerably the resolution, as the the eligibility constraint is a very strong constraint. From the LSSP, the production decision variable $X_{it}$ are directly injected in the second sub-problem: the ASsignment Sub-Problem that will be referred as ASSP (see Figure 7). \\\\}

\begin{figure}[h]
\centering
\resizebox{1\textwidth}{!}{%
\begin{tikzpicture}[node distance=2cm]
    \node (Start) [startstop]{Start};
    \node (Input) [io, right of=Start, xshift= 3cm] {Input Data};
    \node (LSSP) [process, right of=Input, xshift= 3cm] {LSSP};
    \node (SolPool) [io, right of=LSSP, xshift= 4cm] {Solution Pool};
    \node (ASSP) [process, right of=SolPool, xshift= 4cm] {ASSP};
    \node (Dec) [decision, right of=ASSP, xshift= 3cm] {Solution ?};
    \node (Stop) [startstop, right of=Dec, xshift= 3cm] {Stop};
    \draw [arrow] (Start) -- (Input);
    \draw [arrow] (Input) -- (LSSP);
    \draw [arrow] (LSSP) -- (SolPool);
    \draw [arrow] (SolPool) --  node[anchor=south] {Xit} (ASSP);
    \draw [arrow] (ASSP) -- (Dec);
    \draw [arrow] (Dec) -- node[anchor=south] {Yes} (Stop);
    \draw (Dec)-- node[anchor=west] {No} (27,-2,0);
    \draw (27,-2,0)--(16,-2,0);
    \draw [arrow] (16,-2,0)--(SolPool);
   
\end{tikzpicture}
}
\caption{Two-stage Algorithm Flowchart}
\end{figure}
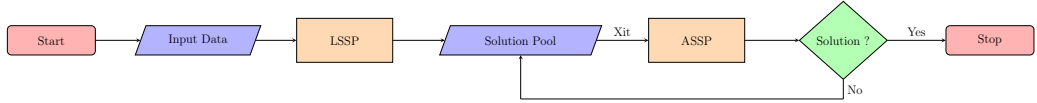
 
\answer{The transformed objective function of the ASSP is to minimize the difference between the production variable $X_{ipt}$ and the production quantities $X_{it}^{LSSP}$ determined by the LSSP. Thanks to this approach, the inventory balance equation and client prioritization modelling are no longer necessary. Thus, the ASSP allows a focus on the machine assignment problem and checks the feasibility of the production plan calculated with the LSSP regarding the eligibility constraint. Below is described the LSSP MIP formulation.}

\subsection{LSSP Model}
\textbf{Indexes and sets}\\
\answer{Indexes and sets used in the MIP (Section \ref{MIP formulation}) remain unchanged.\\}

\textbf{Model parameters}\\
\answer{Model parameters of the MIP model remain unchanged except the following:\\
Matrix $(U_{pt})$ is replaced by $(U_{it})$: the number of molds dedicated to item $i$ in micro-period $t$ to be enforced matrix\\
Matrix $(P_{ipt})$ is no longer necessary.\\}

\textbf{Decision Variables}\\
\answer{Model parameters of the MIP model remain unchanged except the following:\\
$X_{ipt}$ is replaced by $X_{it}$: Quantity of item $i$ to produce in micro-period $t$\\
$Y_{ipt}$ is replaced by $Y_{it}$: Binary variable that equals 1 if there is a new production of an item $i$ at period $t$; 0 otherwise\\
$s_{ipt}$ is replaced by $s_{it}$: Binary variable that equals 1 if there is a mold setup for item $i$ in micro-period $t$; 0 otherwise\\
A new decision variable is also introduced to deal with the mold setup constraint:\\
$\nu_{it}$: Number of molds of item $i$ used in micro-period $t$\\\\
}
The objective function (\ref{2}) remain unchanged.\\

Subject to
\begin{equation}\label{47}
I_{at-1} + \sum_{i=1}^{N_a} X_{it} = I_{at} + \sum_{c \in C} (D_{act} -B_{act} + B_{act-1}) , a = 1..A, t = 1..T
\end{equation}
\begin{equation}\label{48}
I_{at-1}  + \sum_{i=1}^{N_a} X_{it} \ge I_{at} , a = 1..A, t = 1..T
\end{equation}
\begin{equation}\label{49}
I_{at-1} + \sum_{i=1}^{N_a} X_{it} \ge D_{act} - B_{act} + B_{act-1} , a = 1..A, c = C_1 .. C_{n-1}, t = 1..T
\end{equation}
\begin{equation}\label{50}
B_{act} \le D_{act} + B_{act-1} , a = 1..A, c = C_2 .. C_{n}, t = 1..T
\end{equation}\\
The inventory control and customer prioritization constraints (Constraints (6) to (9)) are unchanged compared to the exact method without taking into account the curing press assignment for production quantities.\\
\begin{equation}\label{51}
\nu_{it}\le M Y_{it}, i = 1..N, t = 1..T
\end{equation}
\begin{equation}\label{52}
\nu_{it}\ge Y_{it}, i = 1..N, t = 1..T
\end{equation}
\begin{equation}\label{53}
\sum_{i=1}^{N_a} \nu_{it}\le K_a, a = 1..A, t = 1..T
\end{equation}
\begin{equation}\label{54}
X_{it}=R_{it} \nu_{it}, i = 1..N, t = 1..T
\end{equation}\\
Constraint (\ref{12}) is no longer necessary as there is no curing press in the LSSP. However to anticipate the assignment problem of the ASSP, the production quantities $X_{it}$ are adapted to partially take into account the ignored constraints. Those quantities should be multiples of the daily rate of one mold. Thus, in Constraint (\ref{54}), we are trying to find the best lot sizing with selection of the molds available, without knowing in which curing press these molds will be placed later on. Constraints (\ref{51}) to (\ref{53}) are adapted to match (\ref{54}).\\
\begin{equation}\label{55}
(1-F_{t-1})    s_{it}\ge (1-F_{t-1})   (\nu_{it} - \nu_{it-1}) , i = 1..N, t = 1..T
\end{equation}
\begin{equation}\label{56}
F_{t-1}    s_{it}\ge F_{t-1}   (\nu_{it} - \nu_{it-1-\sum_{o =t-\tau_m}^{t-1}F_{o}}) , i = 1..N, t = 1..T
\end{equation}
\begin{equation}\label{57}
\sum_{i=1}^{N}s_{it}\le \overline{C_t}, t = 1..T
\end{equation}
\begin{equation}\label{58}
\sum_{i=1}^{N}\sum_{t=1}^{T} s_{it}\le \overline{C_h}
\end{equation}\\
The mold setup definition  has to be dealt with differently as we do not consider curing presses in the LSSP. Constraints (\ref{15}) to (\ref{17}) are no longer applicable and replaced by Constraints (\ref{55}) and (\ref{56}). The number of new mold assignment from one period to another is considered as a mold setup. Days-off are also considered in the LSSP with parameter $F_t$.\\
\begin{equation}\label{59}
\sum_{i=1}^{N}\sum_{t=1}^{T} X_{it} \Omega_i\le \sum_{t=1}^{T} V_t + \overline{v_h}
\end{equation}
\begin{equation}\label{60}
\sum_{i=1}^{N}\sum_{t=1}^{T} X_{it} \Omega_i\ge \sum_{t=1}^{T} V_t - \underline{v_h}
\end{equation}
\begin{equation}\label{61}
\sum_{i=1}^{N} X_{it} \Omega_i\le V_t + \overline{v_t}, t = 1..T
\end{equation}
\begin{equation}\label{62}
\sum_{i=1}^{N} X_{it} \Omega_i\ge V_t - \underline{v_t}, t = 1..T
\end{equation}
The weights constraints are unchanged compared to the exact method without taking into account the curing press assignment for production quantities.

\begin{equation}\label{63}
\sum_{i=1}^N Y_{it} \le S, t = 1..T 
\end{equation}

Constraints (\ref{24}) to (\ref{26}) that helps set the $\sigma_{it}$ decision variable in the MIP formulation are no longer useful and simultaneity can be dealt with using Constraint (\ref{63}). The campaign ending modelling remains (Constraints (\ref{27}) and (\ref{28})).\\

\begin{equation}\label{64}
\sum_{i\in N_w}\sum_{t=1}^{T} X_{it} T_i \le \overline{S_w} , w = 1..W
\end{equation}

\begin{equation}\label{65}
\Delta_{idt} \ge \frac{ \nu_{it}}{\varepsilon_{id}} ,i = 1..N, d = 1..N_d, t = 1..T
\end{equation}

\begin{equation}\label{66}
 \frac{ \nu_{it}}{\varepsilon_{id}} \ge 1 - M (1 - \delta_{idt}) ,i = 1..N, d = 1..N_d, t = 1..T
\end{equation}
\begin{equation}\label{67}
 \frac{ \nu_{it}}{\varepsilon_{id}} \le 1 + M \delta_{idt} ,i = 1..N, d = 1..N_d, t = 1..T
\end{equation}

Saturation of upstream assembling workshop is also dealt with in the LSSP without consideration of curing presses in Constraint (\ref{64}). The assembling workshop resources (drums) saturation Constraints (\ref{30}) to (\ref{37}) are taken into account as well and adapted through Constraints (\ref{65}) to (\ref{67}). On the contrary, the minimum production quantity after a mold setup on a given curing press (Constraints (\ref{38}) to (\ref{40})) are not considered in the LSSP.

\begin{equation}\label{68}
\nu_{it} \ge M_{it} , i = 1..N, t = 1..T
\end{equation}
Constraint (\ref{41}) used to integrate quality and industrialisation departments requirements is adapted in Constraint (\ref{68}). However Constraints (\ref{42}) to (\ref{44}) modelling the tire-curing press eligibility and to book a particular curing press are no longer necessary.

\begin{equation}\label{69}
Y_{it}, s_{it}, e_{it},\delta_{idt} \in \{0;1\}
\end{equation}
\begin{equation}\label{70}
X_{it},I_{at},B_{act}, \nu_{it},\Delta_{idt} \in \mathbf{N}
\end{equation}

Finally, Constraint (\ref{69}) defines the boolean decision variables and Constraint (\ref{70}) the production, inventory and number of molds and drums decision variables.

\subsection{ASSP Model}
The ASSP is closer to the MIP formulation described in Section \ref{MIP formulation}. In the ASSP the inventory equation is removed to reduce the model complexity. Instead we use the production quantities $X_{it}^{LSSP}$ from the LSSP as guideline to generate the quantities $X_{ipt}$ to produce in ASSP. Thus the ASSP is only considering the machine assignment decision. It is worth noting that the more restrictive the eligibility matrix ($M_{ap}$), the more complex the ASSP. The objective is then to minimize the gap between the final $X_{ipt}$ quantities and the $X_{it}^{LSSP}$ proposed by the LSSP, which are the ``ideal''  quantities to produce if the the eligibility matrix was an ``all-ones''  matrix. The indexes, sets, model parameters and decision variables (except inventory $I_{at}$ and backorder $B_{act}$ variables that are removed) are the same as in  MIP formulation presented in section \ref{MIP formulation}. The only parameter to be added is the production quantities $X_{it}^{LSSP}$ from the LSSP. The objective function is expressed as follows.\\

Minimize: 
\begin{equation}\label{71}
Z'' = \sum_{i=1}^N \sum_{t=1}^T \big| X_{it}^{LSSP} - \sum_{p=1}^{P} X_{ipt} \big|
\end{equation}

Constraint (\ref{6}) to (\ref{9}) are no longer necessary in ASSP as the demand prioritization has been dealt with in LSSP. The other constraints (Constraints (10) to (46)) are unchanged compared to the MIP formulation presented in Section \ref{MIP formulation}.
\answer{\subsection{Model sizes comparison}}

\answer{The considered problem can be viewed as an extension of the CLSP, as numerous constraints specific to this tire industry case study are added to the classical CLSP. Since CLSP has been proven to be NP-hard by Bitran \cite{bitran1982computational}, our problem is also NP-hard. Besides, to assess the sizes oh the different MIP formulations proposed in this work, a dimensional analysis of the considered problem and model sizes comparison is provided in Table \ref{Dim_analysis}.
}

\begin{table}[h]
\textcolor{black}{
\centering
\begin{tabular}{|c|c|c|c|}
\hline
Model&\# binary v. & \# integer v. & \# constraints\\
\hline
\multirow{2}{*}{Integrated}&\multirow{2}{*}{4.$N.P.T$+2.$N.T$}&\multirow{2}{*}{2.$N.P.T$+3.$N.T$}&21.$N.P.T$+10.$N.T$\\
&&&+2.$P.T$+5.$T$+$W$+4 \\
\hline
\multirow{2}{*}{LSSP}&\multirow{2}{*}{$N.P.T$+6.$N.T$}&\multirow{2}{*}{$N.P.T$+4.$N.T$}&3.$N.P.T$+15.$N.T$\\
&&&+$P.T$+4.$T$+$W$+4 \\
\hline
\multirow{2}{*}{ASSP}&\multirow{2}{*}{4.$N.P.T$+2.$N.T$}&\multirow{2}{*}{2.$N.P.T$}&21.$N.P.T$+3.$N.T$\\
&&&+2.$P.T$+5.$T$+$W$+4 \\
\hline
\end{tabular}
\caption{Model sizes comparison}
\label{Dim_analysis}
}
\end{table}

\answer{
It should be emphasized that both the number of binary and integer variables (in Table \ref{Dim_analysis}, ``\# binary v.'' and ``\# integer v.'' respectively) are reduced in the LSSP in comparison with the Integrated model. The number of constraints (``\# constraints '' in Table \ref{Dim_analysis}) is also considerably shortened. In the ASSP, the number of binary variables is the same as in the Integrated model. There is however less integer variables and cosntraints, and the objective function is easier to solve. Thus, the combination of LSSP and ASSP simplifies the resolution compared to the Integrated model approach.\\}

\answer{In addition, the restrictions on the tire-curing presses eligibility matrix make even more complicated the resolution. Indeed, if we consider a $0-1$ eligibility matrix of size $A.P$, only 32\% of values equals to 1. In other words, each reference of the portfolio can be setup, on average, in less than a third of the total curing presses fleet. It is also to be noted that 17\% of the references can be setup in 10 or less curing presses.}

\section{Computational Results}
\label{Computational Results}
In this section we present an experiment to evaluate the performance of the proposed approaches. We used OPL modelling language and ran the test instances on an Intel Xeon processor with four cores of 3.5 GHz and 64 GB memory. Both MIP and matheuristic were solved using the commercial solver CPLEX 12.9 \cite{cplex2009v12}. Before starting any experiment, the prerequisite is to collect accurate and relevant data, which was a very challenging part of our work. We tested both the exact method and the matheuristic by building a production plan of six weeks as required by the company. However the exact method failed to solve instances with more than 50 items in a reasonable computational time (one hour per production plan of one week allowed). We therefore focus on the results of the matheuristic proposed.
\subsection{Objective function parameter calibration}
\label{OF calib}
To conduct a performance testing campaign on real-world instances, an appropriate parameter calibration for the objective function of the LSSP is necessary, similarly to the work done in \cite{koch2021Capacitated}. First, in order to make fair decision on the five objectives, a normalization has been applied to calculate $\mu_j$ coefficients in the objective function as mentioned in Section 3. Second, the Taguchi design is applied to determine the best values for the coefficients. To do this, performance measures called signal-to-noise (S/N) ratios are calculated. In our case, the term ``signal'' indicates the value that is studied (the mean response) while ``noise'' denotes the uncertainty (standard deviation). The objective function is categorized into three groups by Taguchi method: ``larger is better'', ``smaller is better'' and ``nominal is best''. As we deal with a minimization problem the following ``smaller is better'' S/N ratio is chosen, with $n$ the number of experiments conducted and $Z'$ the  value of the response:
\begin{equation}\label{72}
    S/N =  - 10   log(\frac{\sum_{k=1}^n(Z'_k)^2}{n}) 
\end{equation}

As we deal with three demand classes, five parameters needs to be tuned (backorder of class $C_1$, $C_2$ and $C_3$, overstock and understock). For each parameter, four different values are considered (see Table \ref{Control_factor}). Here two main KPIs can be highlighted: the number of backorders of class $C_1$ and overstock $\overline{S}$. The three others are of a secondary interest for the company. The aggregated weights of the two first KPIs represent 80\% of the total weights, while the remaining 20\% are spread on the other KPIs.\\ 
\begin{table}[h]
\centering
\scalebox{1}{%
\begin{tabular}{|l|c|c|c|c|}
	\hline
	Control Factors&Level 1&Level 2&Level 3&Level 4\\
	\hline
	Backorder $C_1$ ($BC_1$) &20&40&60&80\\
	Overstock ($\overline{S}$)&12&24&36&48\\
	Backorder $C_2$ ($BC_2$)&4&8&12&16\\
	Backorder $C_3$ ($BC_3$) &3&6&9&12\\
	Understock ($\underline{S}$) &1&2&3&4\\
	\hline
\end{tabular}}
\caption{Levels of control factors}
\label{Control_factor}
\end{table}

With five parameters to design and four factor levels $L_{16} (4^5)$ orthogonal array (OA) is chosen to proceed with the experiment. The Design of Experiment is performed with Minitab 19 software \cite{minitab}. Figure 8 represent the signal to noise ratio. Thus the calibration that maximizes the response signal over noise would be the combination of factor levels of Instance 4 (1,4,4,4,4). In Table \ref{Inst_ranking} is presented a ranking of all the sixteen tested instances with an overview over the five measured KPIs. Results of each DoE test instance have been sorted with the following criteria given by the company: i) backorder of class $C_1$; ii) overstock; iii) backorder of class $C_2$; iv) backorder of class $C_3$; v) understock.

\begin{figure}[!h]
	\centering
	\includegraphics[scale=0.6]{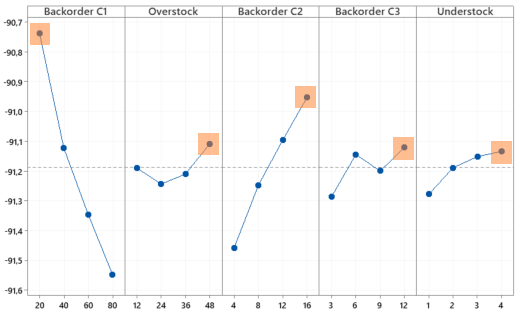}
	\caption{Main effects plots for S/N ratio}
\end{figure}

\begin{table}[h]
\centering
\scalebox{0.7}{%
\begin{tabular}{|l|c|c|c|c|c|c|}
	\hline
	Inst. n°&$BC_1$&$\overline{S}$&$BC_2$&$BC_3$&$\underline{S}$\\
	\hline
	13& \textbf{4(4470)}& 1(130)& 4(47535)& 2(471)& 3(66047)\\
	6& 2(4523)& 2(156)& 1(48585)& 4(464)& 3(66582)\\
	2& 1(4526)& \textbf{2(114)}& 2(44603)& 2(455)& 2(63359)\\
	11& 3(4658)& 3(135)& 1(49973)& 2(469)& 4(67886)\\
	9& 3(4668)& 1(118)& 3(46854)& 4(470)& 2(66012)\\
	5& 2(4823)& 1(132)& 2(46788)& 3(465)& 4(65504)\\
	10& 3(4885)& 2(129)& 4(47314)& 3(471)& 1(66574)\\
	15& 4(4957)& 3(233)& 2(50153)& 4(462)& 1(69084)\\
	14& 4(4974)& 2(158)& 3(49497)& 1(456)& 4(68394)\\
	1& 1(5044)& 1(122)& 1(47105)& 1(473)& 1(65901)\\
	16& 4(5059)& 4(149)& 1(50769)& 3(479)& 2(69830)\\
	8& 2(5102)& 4(132)& 3(45447)& 2(459)& 1(64805)\\
	12& 3(5110)& 4(148)& 2(48256)& 1(470)& 3(67285)\\
	3& 1(5183)& 3(138)& 3(43119)& \textbf{3(451)}& 3(62537)\\
	7& 2(5372)& 3(136)& 4(45358)& 1(455)& 2(64968)\\
	4& 1(5448)& 4(118)& \textbf{4(40813)}& 4(452)& \textbf{4(60528)}\\
	\hline
\end{tabular}}
\caption{Instances ranking regarding company prioritization guidelines }
\label{Inst_ranking}
\end{table}

From this ranking, we can conclude that Instance 13 with coefficient calibration (4,1,4,2,3) from $L_{16}$ Taguchi orthogonal array is the best tuning from the company point of view as it minimizes the number of backorder of class $C_1$. Furthermore, the overstock score of Instance 13 is 130, which is 12\% gap from the best score of overstock (114) provided by Instance 2. It is also to be noted that Instance 4 that minimizes the S/N ratio comes in last position considering the company prioritization guidelines. In the meantime we can notice that Instance 4 minimizes both backorder of class $C_2$ and understock. The relative contribution of each factor to the final result of the objective function is reported in Table \ref{Rel_contrib}.

\begin{table}[h]
\centering
\scalebox{0.7}{%
\begin{tabular}{|l|c|c|c|c|c|c|}
	\hline
	n°&$BC_1$(\%)&$\overline{S}$(\%)&$BC_2$(\%)&$BC_3$(\%)&$\underline{S}$(\%)\\
	\hline
	13&68&10&14&5&3\\
	6&48&29&5&14&4\\
	2&34&40&13&10&3\\
	11&55&33&4&5&4\\
	9&61&12&12&12&2\\
	5&55&16&11&12&5\\
	10&56&23&12&8&1\\
	15&58&26&6&9&1\\
	14&65&20&10&2&3\\
	1&50&30&10&8&3\\
	16&56&34&3&6&1\\
	8&37&45&11&6&1\\
	12&49&39&7&2&2\\
	3&25&45&15&11&4\\
	7&41&37&16&3&2\\
	4&20&48&16&12&4\\
	\hline
\end{tabular}}
\caption{Relative contribution of each factor to the objective function}
\label{Rel_contrib}
\end{table}

\subsection{Case Study and Performance Testing Campaign}
\label{Case Study}
In this section a performance testing campaign is conducted. The benchmark used is build with data provided by the tire manufacturing company. Our hybrid sequential approach is tested on real-world data instances (up to 170 items) on an eight-week planning horizon. The first two weeks of the planning horizon are frozen to stabilize the raw material management process. We therefore tested the matheuristic by building a production plan of six weeks. This production plan is divided in six macro-period of one week (seven working days). For each macro-period a production plan is computed and decision variables are injected from one macro-period to another (see Figure 9). A one-week duration is the smallest time horizon decomposition that we can consider as some constraints such campaign ending or mold setup needs to be dealt with on an entire work-week.\\

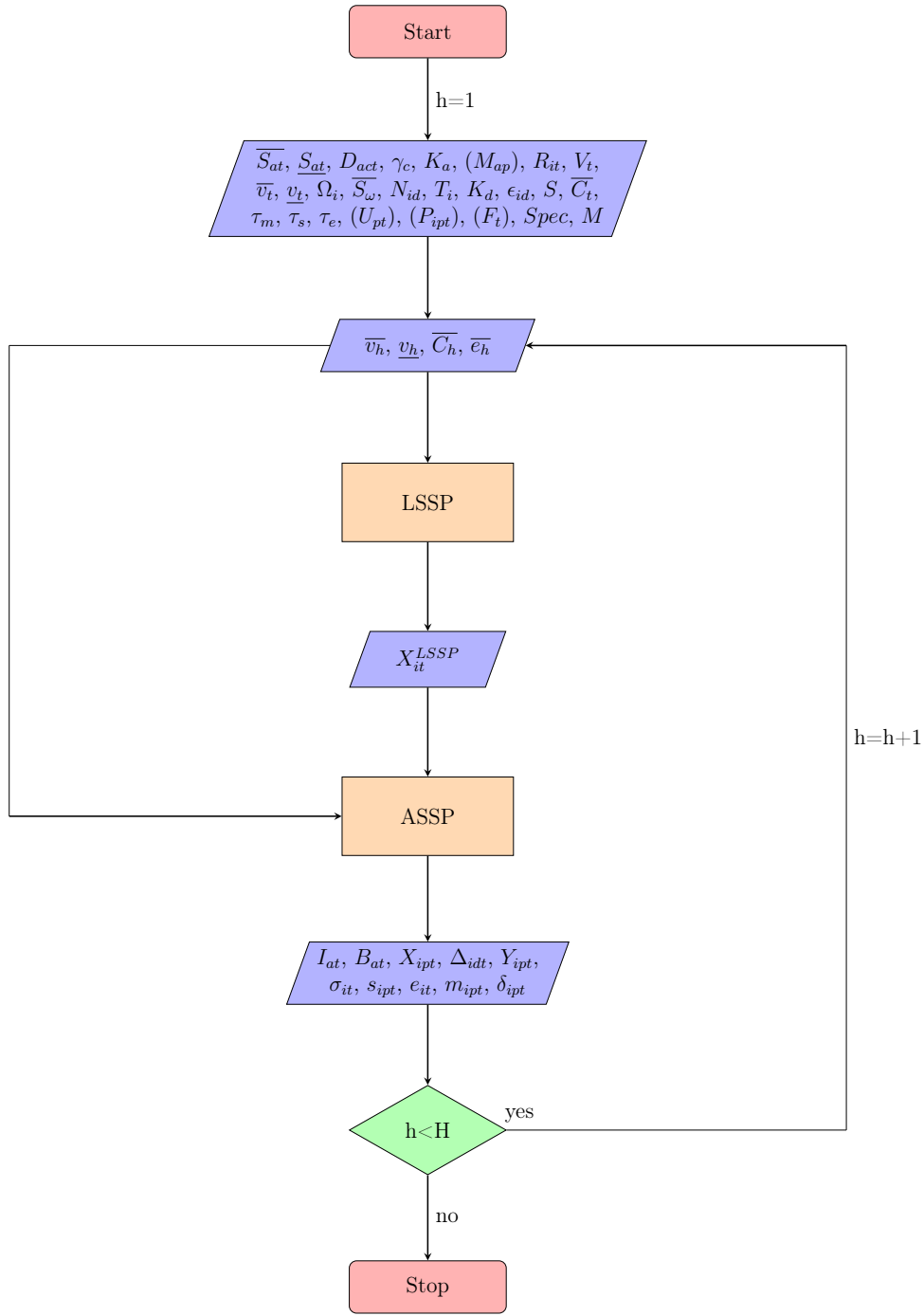
\begin{figure}[!htp]
\centering
\resizebox{0.95\textwidth}{!}{%
\begin{tikzpicture}[node distance=3cm, auto]
    \draw node (start) [startstop]{Start};
    \draw node (in1) [io, below of=start] {
    \begin{varwidth}{17em}
    \begin{center}
    $\overline{S_{at}}$, $\underline{S_{at}}$, $D_{act}$, $\gamma_c$, $K_a$, $(M_{ap})$, $R_{it}$, $V_t$, $\overline{v_t}$, $\underline{v_t}$, $\Omega_i$, $\overline{S_\omega}$, $N_{id}$, $T_i$, $K_d$, $\epsilon_{id}$, $S$, $\overline{C_t}$, $\tau_m$, $\tau_s$, $\tau_e$, $(U_{pt})$, $(P_{ipt})$, $(F_t)$, $Spec$, $M$
    \end{center}
    \end{varwidth}};
    \draw node (in2) [io, below of=in1] {$\overline{v_h}$, $\underline{v_h}$, $\overline{C_h}$, $\overline{e_h}$};
    \draw node (pro2a) [process, below of=in2] {LSSP};
    \draw node (out1) [io, below of=pro2a] {$X_{it}^{LSSP}$};
    \draw node (pro2b) [process, below of=out1] {ASSP};
    \draw node (out2) [io, below of=pro2b] {
    \begin{varwidth}{11em}
    \begin{center}
    $I_{at}$, $B_{at}$, $X_{ipt}$, $\Delta_{idt}$, $Y_{ipt}$, $\sigma_{it}$, $s_{ipt}$, $e_{it}$, $m_{ipt}$, $\delta_{ipt}$ 
    \end{center}
    \end{varwidth}};
    \draw node (dec1) [decision, below of=out2] {h<H};
    \draw node (stop) [startstop, below of=dec1] {Stop};
    \draw [arrow] (start) -- node[anchor=west] {h=1}  (in1) ;
    \draw [arrow] (in1) -- (in2);
    \draw (in2) -- (-8,-6);
    \draw (-8,-6) -- (-8,-15);
    \draw [arrow] (-8,-15) -- (pro2b) ;
    \draw [arrow] (in2) -- (pro2a);
    \draw [arrow] (pro2a) -- (out1);
    \draw [arrow] (out1) -- (pro2b);
    \draw [arrow] (pro2b) -- (out2);
    \draw [arrow] (out2) -- (dec1);
    \draw [arrow] (dec1) -- node[anchor=west] {no} (stop);
    \draw (dec1) -- node[anchor=south, xshift=-3cm] {yes} (8,-21);
    \draw (8,-21) -- node[anchor=west, xshift=0cm] {h=h+1} (8,-6) ;
    \draw [arrow] (8,-6) -- (in2) ;
\end{tikzpicture}
}%
\caption{Decision-support System for Planning and Scheduling Algorithm Flowchart}
\end{figure}

\begin{table}[h]
\textcolor{black}{
    \centering
    \begin{tabular}{|c|c|c|c|c|c|c|c|c|}
    \hline
    Parameter&$S$&$\overline{e_h}$&$\overline{C_t}$&$\overline{C_h}$&$\overline{\tau_m}$&$\overline{\tau_S}$&$\overline{\tau_e}$  \\
    \hline
    Value&43&18&25&5&4&7&4 \\ 
    \hline
    \end{tabular}
    \caption{Flexibility paramaters}
    \label{Table_param_flex}
}
\end{table}

\answer{
In Table \ref{Table_param_flex} is reported an example of flexibility parameters used for this test campaign. These parameters are of a significance importance to measure the flexibility of the factory. The combination of these parameters allow the planning manager to drive properly the inventory turnover of the global portfolio. The simultaneity $S$ and campaign ending $\overline{e_h}$ parameters allow a proper management of material loss and machines saturations in the semi-finite manufacturing workshops. With the mold setup $\overline{C_t}$ and $\overline{C_h}$ parameters, the human ressource capacity of curing workshop saturation is under control. Parameter $\overline{\tau_m}$ ensures production campaigns continuity and a minimum order quantity for each new mold setup to keep setup costs to a minimum.}

 The comparison between the company's solution and the matheuristic is based on the criteria defined in the previous section measured by the company. Namely the number of backorder of class $C_1$ ( $BC_1$), the amount of overstock ($OS$), the number of backorder of class $C_2$ ($BC_2$),  the number of backorder of class $C_3$ ($BC_3$) and finally the value of understock ($US$). However in the day-to-day planning management the company only measures four key performance indicators (KPI):\\
- Two principal KPIs: class $C_1$ backorders ($BC_1$) and overstock ($OS$)\\
- Two secondary KPIs: total number of backorders ($BT$) and understock ($US$)\\
\begin{table}[h!]
\centering
\begin{tabular}{|c|c|c|c|c|c|}
\hline
Dataset   & $OF$ (\%)&  $BC_1$ (\%)& $OS$ (\%)&  $BT$ (\%)& $US$ (\%)\\
\hline
1         & -3   & -15   & -1    & 13    & -4    \\
2         & -36  & -32   & -32   & -38   & -34   \\
3         & -19  & -32   & -23   & -16   & 5     \\
4         & -2   & -34   & -3    & 7     & 0     \\
5         & -30  & -58   & -43   & -6    & -14   \\
6         & -13  & -60   & -7    & -89   & -15   \\
7         & -7   & -54   & 1     & -20   & -9    \\
8         & -8   & -25   & 7     & -40   & -9    \\
9         & -1   & -12   & -21   & 45    & 24    \\
10        & -3   & 1     & -9    & 7     & -4    \\
\hline
\hline
\textbf{Mean}& \textbf{-12} &\textbf{-32} &\textbf{-13} &\textbf{-14} &\textbf{-6} \\
Std. Dev. & 12    & 19    & 15    & 35    & 14    \\
Min.       & -36   & -60   & -43   & -89   & -34   \\
Max.       & -1    & 1     & 7     & 45    & 24    \\
\hline
\end{tabular}
\caption{Comparative table of matheuristic and company solution}
\label{comparison MATH COMP}
\end{table}

\begin{figure}[h]
    \centering
	\includegraphics[scale=0.8]{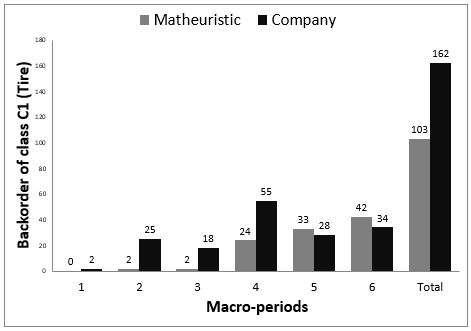}
	\includegraphics[scale=0.8]{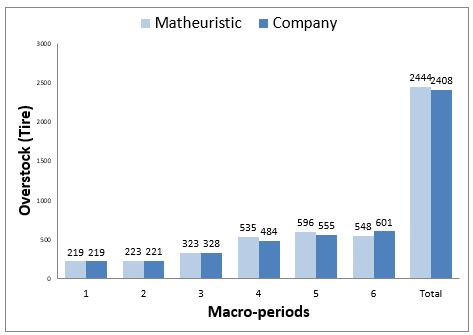}
	\includegraphics[scale=0.8]{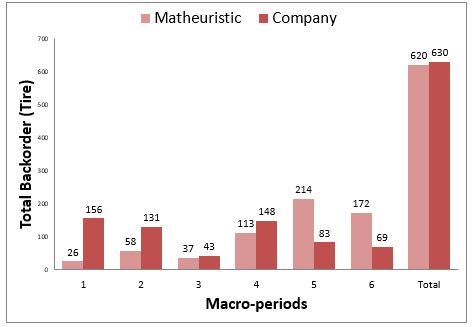}
	\includegraphics[scale=0.8]{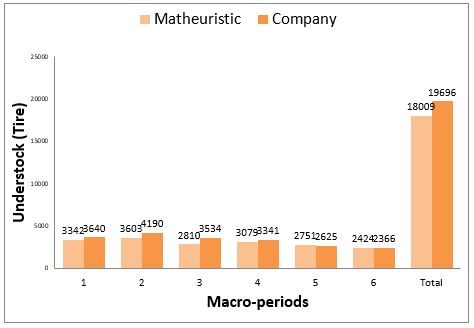}
	\caption{Company KPI for a six-week production plan}
\end{figure}

The results are presented in Table \ref{comparison MATH COMP} in percentage. Ten dataset have been provided by the plant production manager. A negative rate means that the result found with the matheuristic outperforms the company's solution. On the contrary with a positive rate the current solution is better. For instance ``-5'' means that the result of the proposed matheuristic outperfoms the company's solution by five percent. Again the commercial solver CPLEX 12.9 was used to solve the problem, limited to 1800 seconds of computational time for the LSSP and 3600 seconds for the ASSP, with ``MIP emphasis'' parameter set to 1. The first element to be noted is that over ten test campaigns the objective function has always been improved by the matheuristic, with different levels. A significant weight is given to the two principal KPIs, namely the number of backorders of tires of class $C_1$ and the amount of overstock. As a result  $BC_1$ have been reduced by 32\% and overstock by 13\% over the ten datasets, without deterioration of the two secondary KPIs. To give a more complete insight, below is displayed the detailed results of one production plan (see Figure 10). The KPIs are represented macro-period by macro-period and the total amount of tires for each KPI is sumed up in the last bar of the histogram. Thus we can conclude that the decomposition approach presented provides excellent feasible solutions for real-world instances in a reasonable time that outperforms the company's current solution.\\

\subsection{Sensitivity analysis: eligibility constraints}
\answer{As mentionned in the model comparison in Section \ref{Matheuristic}, the eligibility constraints increase the difficukty of the problem. Therefore, a sensitivity analysis of this constraint is proposed. The analysis is conducted on a set of 25 products within the global portfolio. This set of products is a specific range that is strategic for the factory. The idea is to quantify the impact of a new mold technology on the service and inventory level KPIs. This new mold technology implementation would result in a reduction of possibilities in the tire-curing press eligibility matrix. Five different scenarios (see Table \ref{Table_Scenario_Demoul}) are considered and compared to the current situation (Scenario S0). In the current scenario, 30 curing presses are available. In scenarios S1 to S5, 13 to 24 curing presses are made available for the range of products studied. The rest of the eligibility matrix is unchanged for the remaining references of the portfolio.}

\begin{table}[h]
\textcolor{black}{
    \centering
    \begin{tabular}{|c|c|c|c|c|c|c|}
    \hline
    Scenario & S0 & S1 & S2 & S3 & S4 & S5 \\
    \hline
    \# Presses & 30 & 13 & 18 & 20 & 22 & 24 \\
    \hline
    \end{tabular}
    \caption{Scenarios considered for eligibility constraint sensitivity analysis}
    \label{Table_Scenario_Demoul}
    }
\end{table}

\answer{
For each of the 6 scenarios considered, three configurations are tested: low, medium, and high initial inventory level. For each configuration, two values are given: the first one is the absolute value result (in tires) of one particular KPI of one scenario. The second one is the relative value compared to scenario S0. Also, the four KPIs presented in Section \ref{Case Study} are used for each individual scenario. The results of this experiment are presented in Tables \ref{Table_Res_Demoul_1} to \ref{Table_Res_Demoul_6}. The objective function calibration is the same as Section \ref{OF calib}.}

\begin{table}\footnotesize
\textcolor{black}{
\centering
\begin{tabular}{|c|c|c|c|c|c|}
\hline
\multirow{2}{*}{Scenario \#} & \multirow{2}{*}{KPI} & Configuration 1 &Configuration 2 & Configuration 3  \\
\cline{3-5}
&&Tires&Tires&Tires\\
\hline
\multirow{4}{*}{Scenario S0}  &  $BC_1$ & 2563 & 226 &294\\
&  $OS$ & 0 & 543 & 902\\
 &  $BT$ & 43103 & 6298 &3758\\
 &  $US$ & 51849 &30676 &6215\\
\hline
\end{tabular}
    \caption{Reference scenario S0 results}
    \label{Table_Res_Demoul_1}
    }
\end{table}
 
\answer{In Configuration 1 (low initial inventory level) for scenario S0, the backorder of class $C_1$, the total backorder $BT$ and the understock $US$ levels are very high while the overstock $OS$ is null (See Table \ref{Table_Res_Demoul_1}). In Configuration 2 and 3 (medium and high initial inventory level), the total backorder $BT$ and the understock $US$ levels improve in coherence with the initial inventory. Besides, with increasing inventory level from Configuration 1 to 3, overstock $OS$ level arise consequently. It is to be noted that backorder of class $C1$ are contained respectively to 226 and 294 in Configuration 2 and 3.}

\begin{table}[h]\footnotesize
\textcolor{black}{
\centering
\begin{tabular}{|c|c|c|c|c|c|c|c|}
\hline
\multirow{2}{*}{Scenario \#} & \multirow{2}{*}{KPI} & \multicolumn{2}{c|}{Configuration 1} &\multicolumn{2}{c|}{Configuration 2} & \multicolumn{2}{c|}{Configuration 3}   \\
\cline{3-8}
&&Tires&Gap S0&Tires&Gap S0&Tires&Gap S0\\
\hline
\multirow{4}{*}{Scenario S1}  &  $BC_1$ & 5209 & 103\% & 840 & 272\% & 912 & 210\% \\
&  $OS$ & 16456 & NA & 1561 & 187\% & 5618 & 523\% \\ 
 &  $BT$ & 52607 & 22\% & 7516 & 19\% & 6779 & 80\% \\ 
&  $US$ & 61980 & 6\% & 34145 & 11\% & 11016 & 77\% \\
\hline
\end{tabular}
    \caption{Scenario S1 results}
    \label{Table_Res_Demoul_2}
    }
\end{table}

\answer{In Table \ref{Table_Res_Demoul_2} are presented the results of Scenario S1, that represent a significative flexibility degradation as only 13 curing presses are available for the range of products studied within the global portfolio. In every configurations, both the backorder of class $C_1$ and the overstock $OS$ levels increase dramatically. In particular, the backorder of class $C_1$ level more than triple in Configuration 3 whereas initial inventory level are high. Similarly, even with low initial inventory level in Configuration 1, the overstock $OS$ level soars, increasing from 0 to 16456 tires from Scenario S0 to Scenario S1. Here, the notation `` NA '' indicates the impossible division by zero. Moreover, it should be noted that the total backorder $BT$ and the understock $US$ are contained to an increase of less than 25\% in Configuration 1 and 2 while rising to around 80\% in Configuration 3.}

\begin{table}[h]\footnotesize
\textcolor{black}{
    \centering
\begin{tabular}{|c|c|c|c|c|c|c|c|}
\hline
\multirow{2}{*}{Scenario \#} & \multirow{2}{*}{KPI} & \multicolumn{2}{c|}{Configuration 1} &\multicolumn{2}{c|}{Configuration 2} & \multicolumn{2}{c|}{Configuration 3}   \\
\cline{3-8}
&&Tires&Gap S0&Tires&Gap S0&Tires&Gap S0\\
\hline
 \multirow{4}{*}{Scenario S2}  &  $BC_1$                   & 3681             & 44\%              & 486                & 115\%              & 793               & 170\%               \\  
                     &  $OS$                   & 967              & NA                & 1286               & 137\%              & 3291              & 265\%              \\  
                     &  $BT$                   & 45611            & 6\%               & 7065               & 12\%               & 5559              & 48\%                \\  
                     &  $US$                   & 54933            & 6\%               & 32214              & 5\%                & 8910              & 43\%                \\ \hline                    
\end{tabular}
    \caption{Scenario S2 results}
    \label{Table_Res_Demoul_3}
    }
\end{table}

\begin{table}[h]\footnotesize
\textcolor{black}{
    \centering
\begin{tabular}{|c|c|c|c|c|c|c|c|}
\hline
\multirow{2}{*}{Scenario \#} & \multirow{2}{*}{KPI} & \multicolumn{2}{c|}{Configuration 1} &\multicolumn{2}{c|}{Configuration 2} & \multicolumn{2}{c|}{Configuration 3}   \\
\cline{3-8}
&&Tires&Gap S0&Tires&Gap S0&Tires&Gap S0\\
\hline
\multirow{4}{*}{Scenario S3}  &  $BC_1$                   & 3787             & 48\%              & 347                & 54\%               & 785               & 167\%                \\  
                     &  $OS$                   & 866              & NA                & 567                & 4\%                & 2999              & 232\%             \\  
                     &  $BT$                   & 45999            & 7\%               & 6647               & 6\%                & 5067              & 35\%                  \\  
                     &  $US$                   & 54933            & 6\%               & 31768              & 4\%                & 8095              & 30\%                \\ \hline
\end{tabular}
    \caption{Scenario S3 results}
    \label{Table_Res_Demoul_4}
    }
\end{table}

\answer{In Table \ref{Table_Res_Demoul_3} and \ref{Table_Res_Demoul_4} are presented the results of Scenarios S2 and S3, in which the number of curing presses available are respectively 18 and 20. It should be emphasized that the trend of Scenario S1 is the same for each KPI, though with a smaller scale.}

\begin{table}[h]\footnotesize
\textcolor{black}{
    \centering
\begin{tabular}{|c|c|c|c|c|c|c|c|}
\hline
\multirow{2}{*}{Scenario \#} & \multirow{2}{*}{KPI} & \multicolumn{2}{c|}{Configuration 1} &\multicolumn{2}{c|}{Configuration 2} & \multicolumn{2}{c|}{Configuration 3}   \\
\cline{3-8}
&&Tires&Gap S0&Tires&Gap S0&Tires&Gap S0\\
\hline
\multirow{4}{*}{Scenario S4}  &  $BC_1$                   & 2272             & \textbf{-11\%}& 323                & 43\%               & 271               & \textbf{-8\%}\\  
                     &  $OS$                   & 36               & NA                & 653                & 20\%               & 1149              & 27\%                 \\  
                     &  $BT$                   & 43282            & 0\%               & 6214               & \textbf{-1\%}& 3482              & \textbf{-7\%} \\  
                     &  $US$                   & 51341            & \textbf{-1\%}& 29490              & \textbf{-4\%}& 6070              & \textbf{-2\%}\\ \hline
\end{tabular}
    \caption{Scenario S4 results}
    \label{Table_Res_Demoul_5}
    }
\end{table}

\answer{Table \ref{Table_Res_Demoul_5} shows the results of Scenario S4 test instances that allow 22 curing presses. It should be highlighted that the backorder of class $C_1$ level derceased by around 10\% in both Configurations 1 and 3 while still arise to 43\% in Configuration 2. The overstock $OS$ level logically increases from Configuration 1 to 3 in coherence with initial inventory level and is notably contained to an increase of 27\% in Configuration 3. Please also note that the total backorder $BT$ and the understock $US$ levels are slightly improved in a range varying from 0\% to 7\% in every configuration.}

\begin{table}[h]\footnotesize
\textcolor{black}{
    \centering
\begin{tabular}{|c|c|c|c|c|c|c|c|}
\hline
\multirow{2}{*}{Scenario \#} & \multirow{2}{*}{KPI} & \multicolumn{2}{c|}{Configuration 1} &\multicolumn{2}{c|}{Configuration 2} & \multicolumn{2}{c|}{Configuration 3}   \\
\cline{3-8}
&&Tires&Gap S0&Tires&Gap S0&Tires&Gap S0\\
\hline
\multirow{4}{*}{Scenario S5}  &  $BC_1$                   & 2312             & \textbf{-10\%}& 210 & \textbf{-7\%}& 275& \textbf{-6\%}\\  
&  $OS$& 16 & NA& 528& \textbf{-3\%} & 1026 & 14\% \\  
 &  $BT$& 42335 & \textbf{-2\%} & 6145 & \textbf{-2\%}& 3992& 6\%\\  
 &  $US$& 50954& \textbf{-2\%} & 29179& \textbf{-5\%}& 6540& 5\%\\ \hline
\end{tabular}
    \caption{Scenario S5 results}
    \label{Table_Res_Demoul_6}
    }
\end{table}
 
\answer{Lastly, in Table \ref{Table_Res_Demoul_6} are reported the results of Scenario S5 in which 24 curing presses are made available. It is to be noted that the backorder of class $C_1$ level is limited in a range of 6\% to 10\% for every configuration compared to Scenario S0. The overstock $OS$ level is lightly improved in Configuration 2 and contained to a rise of 14\% in Configuration 3. Furthermore, the total backorder $BT$ and the understock $US$ level are improved in a range from 2\% to 5\% in both Configurations 1 and 2 while limited to less than 6\% in Configuration 3.}

\answer{To put it in a nutshell, the main driver for the sensitivity analysis is to consider the backorder of class $C_1$ and the overstock $OS$ KPIs. From this, we can conclude that only Scenarios S4 and S5 with respectively 22 and 24 machines available would be acceptable, although only Scenario S4 performs as well as the current situation of Scenario S0. This conclusion is confirmed by the analysis of the total backorder $BT$ and the understock $US$ results. In other words, with 20 or less machines available for this range of products, in combination with the rest of the portfolio that compete for the same ressources, all the KPIs explode and the available capacity is tightened, whichdoes not allow to achieve satisfactory service and inventory level. It should also be emphasized that in Scenario S1, S2 and S3, even with low initial inventory, the $OS$ score increase dramatically to satisfy both tonnage and eligibility constraints.}

\section{Conclusion}
\label{conclusion}
\answer{
In this paper, a dedicated lot sizing problem inspired from the off-the-road tire industry is addressed. This problem is a single-level multi-item multi-machine simultaneous lot sizing and scheduling problem with dynamic demand, backlogging, sequence independent setup times and specific constraints. A mathematical model, which describes this complex industrial case and presents new constraints and an original sequential approach is proposed. To the best of our knowledge, modelling a client prioritization, upstream workshop resources and the number of campaign ending have never been dealt with in the literature simultaneously. This allows to solve a real-world problem with up to 170 items on 70 machines and 42 periods. The contribution of this work is threefold.
First, an industrial case is tackled and the efficiency of the proposed method to solve it is demonstrated.
Second, a new MIP formulation of a very complex problem in the off-the-road tire industry is provided.
Finally, the proposed hybrid sequential approach managed to solve the real-world problem in reasonable computational time and outperformed the current solution of the company. The two most important KPIs for the management have been optimized of respectively 32\% for the backorders of class $C_1$ and 13\% for the overstock. Moreover, the resolution time have been reduced significantly.}

\answer{
Besides, this study drives two main outputs for the factory management team. First, the proposed framework is core to a decision support system for planning and scheduling that both model the company’s production system and generate production plans. This research therefore allows the company to automate their production planning process, analyse data and drive continuous improvement. Second, the proposed framework has been and can be used as a management support and performance evaluation tool to the head of the factory to help them take better decisions. The sensitivity analysis provided in the result section is one example of numerous applications that can be conducted using this framework. The plant management is now able to elaborate quantitative scenarios about their future strategic decisions to invest in and modernize their industrial system. Also, the very promising results yielded by the matheuristic leads the company to study the implementation of our resolution method with a dedicated software.}

\answer{
Finally, one strong stance of the study is to solve a very complex case study using mathematical modelling instead of classical more widespread metaheuristic procedure. Therefore, comparing those two approaches based on the literature review may be a very insightful future research perspective. In addition, the close relationship between the assembling and curing workshops has been highlighted in Section \ref{Problem description}. Approaching the whole problem using a hybrid flow shop environment with two stages could be a very enlightening outlook. Besides, the specific constraints introduced in this work contribute to the extension of the already considerable lot sizing types and variants. As these constraints are directly extracted from a real-world problem, we believe that future works can consider different combinations that may arise in other companies and industries.}

\bibliography{biblio}

\end{document}